\documentclass[11pt]{article}

\usepackage{amsfonts}
\usepackage{amssymb} 
\usepackage{german}

\def\II{I\hspace{-0.1cm}I}

\def\negthickspace{\!\!\!}

\newcommand{\nicefrac}[2]
{\leavevmode \kern.1em\raise.5ex\hbox{\the\scriptfont0 #1}
             \kern-.1em/\kern-.15em\lower.25ex
             \hbox{\the\scriptfont0 #2}}

\newtheorem{Theo}{Theorem}{\alph{enumi}}
\newenvironment{theorem}{\begin{Theo}\hspace{-0.2cm}: }{\end{Theo}}

\newtheorem{Pro}{Proposition}{\alph{enumi}}
\newenvironment{proposition}{\begin{Pro}\hspace{-0.2cm}: }{\end{Pro}}

\newtheorem{Co}{Corollary}{\alph{enumi}}
\newenvironment{corollary}{\begin{Co}\hspace{-0.2cm}: }{\end{Co}}

\newtheorem{As}{Assumption}{\alph{enumi}}

\newtheorem{Le}{Lemma}{\alph{enumi}}

\newtheorem{Fo}{Folgerung}{\alph{enumi}}

\newtheorem{De}{Definition}{\alph{enumi}}
\newenvironment{definition}{\begin{De}\hspace{-0.2cm}:\rm }{\end{De}}

\newtheorem{Be}{Bemerkung}{\alph{enumi}}

\newtheorem{Ex}{Example}{\alph{enumi}}
\newenvironment{example}{\begin{Ex}\hspace{-0.2cm}:\rm }{\end{Ex}}

\begin{document} 

\begin{center}
{\Large{\sc On two-dimensional immersions \\[0.1cm] of prescribed mean curvature in $\mathbb R^n$}}\\[1.2cm]
{\large{\sc Matthias Bergner, Steffen Fr\"ohlich}}\\[1cm]
{\small\bf Abstract}\\[0.4cm]
\begin{minipage}[c][2cm][l]{12cm}
{\small We consider two-dimensional immersions of disc-type in $\mathbb R^n.$ We focus well known classical concepts and study the nonlinear elliptic systems of such mappings. Using an Osserman-type condition we give a priori-estimates of the principle curvatures for certain graphs in $\mathbb R^4$ with prescribed mean curvature.}
\end{minipage}
\end{center}
{\bf\large Introduction}\\[0.3cm]
In 1917, S. Bernstein proved the following theorem:\\[0.1cm]
{\it If the twice continuously differentiable function $f:\mathbb R^2\to\mathbb R$ represents a graph of a complete minimal surface in $\mathbb R^3,$ then $f$ is linear, i.e. the minimal surface is a plane.}\\[0.1cm]
This result indicates the difference between linear partial differential equations and the non-linear minimal surface equation. Over the years much work has been devoted in finding new proofs and generalizations.\\[0.1cm]
One method which goes back to works of E. Heinz and F. Sauvigny will be used in the paper at hand: Using isothermal parameters, a minimal surface $X$ and its unit normal mapping $N$ in $\mathbb R^3$ satisfy the nonlinear elliptic systems
  $$\triangle X=0,\quad \triangle N=-2KWN.$$
Here, $K$ denotes the Gaussian curvature and $W$ the surface area element of the immersion $X.$ If one can establish certain bounds on the size of $X$ and $N,$ the gradient estimates in \cite{Heinz_01} ensure upper bounds for the first and second derivatives of $X$ in terms of given a priori-data. A Harnack-type inequality for pseudo-holomorphic functions gives lower bounds for $W.$ Altogether, this gives estimates for the principle curvatures of the immersion, namely in such a way that if the surface growths to become complete, the principle curvatures vanish identically.\\[0.1cm]
This method generalizes to various variational problems if the differential systems satisfy a structure condition of the form
  $$|\triangle Z|\le\mbox{const}\,|\nabla Z|^2\,.$$
For example, the mean curvature system $\triangle X=2HWN$ and the above system for $N$ satisfy this assumption. Some problems arise in higher codimension.\\[0.1cm]
In 1964, Osserman \cite{Osserman_01} proved that a complete two-dimensional minimal surface in $\mathbb R^n$ is a plane if all of its normal vectors make a certain positive angle with a fixed axis in space. The assumption that {\it all normals make an acute angle with some fixed direction} makes the difference to the Bernstein theorem in $\mathbb R^3$ where it is sufficient to assume that the normal of the immersion omits a certain neighbourhood of some point of the unit sphere.\\[0.1cm]
Osserman's method is based essentially on results of complex analysis due to the fact that the components of a minimal surface are harmonic functions.
\newpage\noindent
Although it is rather difficult to realize Osserman's condition, it is an important fact that any condition of this type is necessary for curvature estimates.\\[0.1cm]
In this paper we want to examine how to apply the methods developed in \cite{Heinz_01} and \cite{Sauvigny_01} to establish certain a priori-estimates for the principle curvatures of graphs with prescribed mean curvature. In the first part, we concentrate on some geometric basics of two-dimensional immersion in $\mathbb R^n$ which are motivated from \cite{Brauner_01}, and we will derive nonlinear elliptic differential systems to describe the surfaces analytically. The second part contains the mentioned estimates.
\subsection{Surfaces in $\mathbb R^n$}
The essential principles of the following are motivated from \cite{Blaschke_Leichtweiss_01}, \cite{Brauner_01}, \cite{Osserman_01}, and \cite{Weyl_01}.\\[0.1cm]
Let
  $$B:=\{(u,v)\in\mathbb R^2\,:\,u^2+v^2<1\}$$
denote the open unit disc, $\overline B\subset\mathbb R^2$ its topological closure. For integer $n\ge 3$ we consider immersions
  $$\begin{array}{l}
      X=X(u,v)\in C^{2+\alpha}(B,\mathbb R^n)\cap C^0(\overline B,\mathbb R^n),\quad\alpha\in(0,1), \\[0.1cm]
      X(u,v)=(x^1(u,v),x^2(u,v),\ldots,x^n(u,v)),
    \end{array}$$
with the property
  $$\mbox{rank}\,\partial X(u,v)
    \equiv\mbox{rank}
      \left(
        \begin{array}{cc}
          x_u^1(u,v) & x_v^1(u,v) \\[0.1cm] \vdots & \vdots \\[0.1cm] x_u^n(u,v) & x_v^n(u,v)
        \end{array}
      \right)
    =2\quad\mbox{for all}\ (u,v)\in B.$$
This means that $X_u=X_u(u,v)$ and $X_v=X_v(u,v)$ are linearly independent at every point $(u,v)\in B$ and, therefore, span the two-dimensional tangent plane
  $${\mathcal T}_X(w):=\mbox{Span}\,\{X_u(w),X_v(w)\}\quad\mbox{at}\ w\in B.$$
There exists a normal space ${\mathcal N}_X(w)$ spanned by $n-2$ linearly independent vectors $N_1,N_2,\ldots,N_{n-2}$ which can be assumed to satisfy the orthonormality relations
  $$N_\Sigma\cdot N_\Theta^t
    =\delta_{\Sigma\Theta}
    :=\left\{
        \begin{array}{l}
          1\quad\mbox{if}\ \Sigma=\Theta \\[0.1cm]
          0\quad\mbox{if}\ \Sigma\not=\Theta
        \end{array}
      \right.\quad\mbox{for all}\ \Sigma,\Theta\in\{1,2,\ldots,n-2\}.$$
Here, by the upper $t$ we denote the transposed vector.\\[0.1cm]
The choice of the normal vectors $N_\Sigma(w),$ $\Sigma=1,\ldots,n-2,$ is not unique.
\begin{example}
The unit vectors
  $$\begin{array}{lll}
      N_1\negthickspace
      & := & \negthickspace\displaystyle
             \frac{1}{\sqrt{1+|\nabla\varphi_1|^2}}\,(-\varphi_{1,x},-\varphi_{1,y},1,0,\ldots,0), \\[0.6cm]
      N_2\negthickspace
      & := & \negthickspace\displaystyle
             \frac{1}{\sqrt{1+|\nabla\varphi_2|^2}}\,(-\varphi_{2,x},-\varphi_{2,y},0,1,\ldots,0),\ \ldots
    \end{array}$$
are normal to the graph $X(x,y)=(x,y,\varphi_1(x,y),\ldots,\varphi_n(x,y)).$
\end{example}
{\bf Assumption: }{\it There exists an orthonormal basis $N_1(w),N_2(w),\ldots,N_{n-2}(w)$ of the normal space
$\mathcal N_X(w)$ such that $N_\Sigma\in C^{1+\alpha}(B,S^{n-1}).$} \\[0.3cm]
In the following we will use such regular bases.
\subsubsection{Differential geometric foundations}
We consider surfaces $X=X(u,v)$ immersed into the Euclidean space $\mathbb R^n.$ Inserting $X=X(u,v)$ into the element $ds^2:=(dx^1)^2+\ldots+(dx^n)^2$ yields
  $$ds^2=h_{11}\,du^2+2h_{12}\,dudv+h_{22}\,dv^2=h_{ij}\,du^idu^j$$
(note the summation convention in the last term) with the first fundamental form
    $$I:=(h_{ij})_{i,j=1,2}\equiv(X_{u^i}\cdot X_{u^j}^t)_{i,j=1,2}\subset\mathbb R^{2\times 2}\,.$$
Using conformal parameters $(u,v)\in B,$ the first fundamental form appears as
  $$h_{11}=W=h_{22}\,,\quad
    h_{12}=0
    \quad\mbox{in}\ B$$
with the surface area element $W:=\sqrt{h_{11}h_{22}-h_{12}^2}.$\\[0.1cm]
The second fundamental form $\II_\Sigma=(L_{\Sigma,ij})_{i,j=1,2}$ w.r.t. the unit normal $N_\Sigma$ is defined as
  $$L_{\Sigma,ij}:=-X_{u^i}\cdot N_{\Sigma,u^j}^t=X_{u^iu^j}\cdot N_\Sigma^t\,,\quad i,j=1,2.$$
Consider the form
  $$(L_{\Sigma,i}^k)_{i,k=1,2}:=(L_{\Sigma,ij}h^{jk})_{i,k=1,2}\in\mathbb R^{2\times 2}$$
with the coefficients $h^{ij}$ of the inverse of the first fundamental form, i.e. $h_{ij}h^{jk}=\delta_i^k$ with the Kronecker symbol $\delta_i^k.$
\begin{definition}
The {\it mean curvature} $H_\Sigma$ {\it in direction} $N_\Sigma$ is defined as
  $$H_\Sigma:=\frac{1}{2}\,\mbox{trace}\,(L_{\Sigma,i}^k)_{i,k=1,2}
             =\frac{L_{\Sigma,11}h_{11}-2L_{\Sigma,12}h_{12}+L_{\Sigma,22}h_{22}}{2(h_{11}h_{22}-h_{12}^2)}\,.$$
\end{definition}
Besides the mean curvature we introduce the {\it Gaussian curvature w.r.t. to} $N_\Sigma$ as
  $$K_\Sigma
    :=\mbox{det}\,(L_{\Sigma,i}^k)_{i,k=1,2}
    =\frac{L_{\Sigma,11}L_{\Sigma,22}-L_{\Sigma,12}^2}{h_{11}h_{22}-h_{12}^2}\,.$$
The {\it principle curvatures} $\kappa_{\Sigma,1},$ $\kappa_{\Sigma,2}$ w.r.t. $N_\Sigma$ are the eigenvalues of the form $(L_{\Sigma,i}^k)_{i,k=1,2},$ that is
  $$H_\Sigma=\frac{\kappa_{\Sigma,1}+\kappa_{\Sigma,2}}{2}\,,\quad
    K_\Sigma=\kappa_{\Sigma,1}\kappa_{\Sigma,2}\,.$$
\subsubsection{The differential equations}
Let $N_\Sigma,$ $\Sigma=1,\ldots,n-2,$ be an orthonormal basis of the normal space. We rewrite the derivatives $N_{\Sigma,u^i}$  and $X_{u^iu^j}$ for and $i,j=1,2$ in terms of the moving frame $X_u,X_v,N_1,\ldots,N_{n-2}.$ The following is discussed more detailed in \cite{Brauner_01}.
\begin{proposition}
Let $N_\Sigma,$ $\Sigma=1,\ldots,n-2,$ be an orthonormal basis of the normal space. Then
  $$N_{\Sigma,u^i}=-L_{\Sigma,ij}h^{jk}X_{u^k}^t+\sigma_{\Sigma,i}^\Theta N_\Theta\,,
    \quad i=1,2,\ \Sigma=1,2,\ldots,n-2,$$
with the {\it torsion coefficients}
  $$\sigma_{\Sigma,i}^\Theta
    :=\left\{
        \begin{array}{cl}
          N_{\Sigma,u^i}\cdot N_\Theta^t & \mbox{if}\ \Sigma\not=\Theta \\[0.1cm]
          0 & \mbox{if}\ \Sigma=\Theta
        \end{array}
      \right.,$$
as well as
 $$X_{u^iu^j}=\Gamma_{ij}^kX_{u^k}+\sum_{\Sigma=1}^{n-2}L_{\Sigma,ij}N_\Sigma$$
with the Christoffel symbols $\Gamma_{ij}^k:=\frac{1}{2}\,h^{kl}(h_{jl,i}+h_{li,j}-h_{ij,l}),$ $h_{ij,k}:=h_{ij,u^k}\,.$
\end{proposition}
\subsubsection{An elliptic system for $X$}
Let $(u,v)\in B$ be conformal parameters. Then, there hold $\Gamma_{11}^1+\Gamma_{22}^1=0,$ $\Gamma_{11}^2+\Gamma_{22}^2=0.$ Thus, the Gauss equations from Proposition 1 yield
  $$\triangle X
    =(\Gamma_{11}^1+\Gamma_{22}^1)X_u
     +(\Gamma_{11}^2+\Gamma_{22}^2)X_v
     +\sum_{\Sigma=1}^{n-2}(L_{\Sigma,11}+L_{\Sigma,22})N_\Sigma
    =2\sum_{\Sigma=1}^{n-2}H_\Sigma WN_\Sigma\,.$$
This system is called the {\it mean curvature system.}\\[0.1cm]
Let $|H_1|+|H_2|+\ldots+|H_{n-2}|\le h_0$ with a real constant $h_0\in(0,+\infty).$ We estimate
  $$|\triangle X|
    \le 2h_0|X_u||X_v|
    \le h_0|\nabla X|^2$$
due to $W=\sqrt{h_{11}h_{22}}=|X_u||X_v|.$ Note the quadratic growth of the gradient on the right hand side. Elliptic systems of this structure are the matter of investigations in \cite{Heinz_01}.
\subsection{A priori-estimates for graphs in $\mathbb R^4$}
In this part we consider immersions
  $$\begin{array}{l}
      X=X(u,v)\in C^{2+\alpha}(B,\mathbb R^4)\cap C^0(\overline B,\mathbb R^4),\quad\alpha\in(0,1), \\[0.1cm]
      X(u,v)=(x^1(u,v),\ldots,x^4(u,v)).
    \end{array}$$
Let a function $\overline \mathcal H(X):\mathbb R^4 \longrightarrow \mathbb R^4$ be given. We define a scalar function ${\mathcal H}(X,Z):\mathbb R^4\times S^3\to\mathbb R$ by
  $${\mathcal H}(X,Z):=\overline \mathcal H(X)\cdot Z^t\quad\mbox{for}\ X\in\mathbb R^4\,,\ Z\in S^3\,.$$
We now consider surfaces of prescribed mean curvature ${\mathcal H}={\mathcal H}(X,Z)$ in the following sense.
\begin{definition}
A surface $X=X(u,v)$ is called a (conformally parametrized) surface of prescribed mean curvature ${\mathcal H}\,:\,\mathbb R^4\times S^3\to\mathbb R$ if it satisfies the system
  $$\begin{array}{l}
      \triangle X=2\mathcal H(X,N_1)W N_1+2\mathcal H(X,N_2)W N_2\,,  \\[0.1cm]
      |X_u|^2=W=|X_v|^2\,,\quad X_u\cdot X_v^t=0\qquad\mbox{in}\ B
    \end{array}$$
for an orthonormal basis $N_1,N_2$ of the normal space.
\end{definition}
{\it Remarks:}
\begin{itemize}
\item[1.]
By the definition of $\mathcal H$, the above system is satisfied for all orthonormal bases of the normal space if it is satisfied for at least one orthonormal basis.
\item[2.]
Using conformal parameters the geometric mean curvature $H_\Sigma(X)\equiv H(X,N_\Sigma)$ of $X$ in direction of the unit normal $N_\Sigma$ can be calculated by
  $$H_\Sigma(X)
    =\frac{(h_{11} X_{uu}-2h_{12}X_{uv}+h_{22}X_{vv})\cdot N_\Sigma^t}{2(h_{11}h_{22}-h_{12}^2)}
    =\frac{\triangle X\cdot N_\Sigma^t}{2W}\,.$$
Thus, by section 1.3, the geometric mean curvature $H_\Sigma(X)$ is equal to $\mathcal H(X,N_\Sigma)$ if the normal basis is continuously differentiable.
\end{itemize}
Our main result is the following.
\begin{theorem}\\[0.1cm]
{\bf Assumptions: }Let $X\in C^{2+\alpha}(B,\mathbb R^4)\cap C^0(\overline B,\mathbb R^4),$ $\alpha\in(0,1),$ be a surface of prescribed mean curvature ${\mathcal H}(X,Z)\in C^0(\mathbb R^4\times S^3,\mathbb R)$ with the following properties:
\begin{itemize}
\item[(A1)]
The immersion $X=X(u,v)$ is a positively oriented conformal reparametrization of a graph of the form
  $$X(x,y)=(x,y,\varphi(x,y),\psi(x,y)),\quad(x,y)\in\overline B_R,$$
where $B_R:=\{(x,y)\in R^2\,:\,x^2+y^2<R^2\},$ and $\varphi,\psi\in C^{2+\alpha}(\overline B_R,\mathbb R).$
\item[(A2)]
The mean curvature ${\mathcal H}={\mathcal H}(X,Z)$ satisfies
  $$|{\mathcal H}(X,Z)|\le\frac{h_0}{2}\quad\mbox{for all}\ X\in\mathbb R^4\,,\ Z\in S^3\,,$$
as well as
  $$|\mathcal H(X_1,Z_1)-\mathcal H(X_2,Z_2)|\leq h_1|X_1-X_2|^\alpha+h_2|Z_1-Z_2|
    \quad\mbox{for all}\ X_1,X_2\in\mathbb R^4\,,\ Z_1,Z_2\in S^3$$
with constants $h_0,h_1,h_2\in [0,+\infty).$
\item[(A3)]
The surface represents a geodesic disc of radius $r>0$ and with center $X(0,0)=(0,0,0,0).$
\item[(A4)]
With a real constant $d_0>0,$ the Dirichlet energy can be estimated by
  $${\mathcal D}[X]
    :=\int\hspace{-0.25cm}\int\limits_{\hspace{-0.3cm}B}|\nabla X(u,v)|^2\,dudv
    \le d_0r^2\,.$$
\item[(A5)]
At every point $w\in B,$ each normal vector of the immersion makes an angle of at least $\omega>0$ 
with the $x_1$-axis.
\end{itemize}
{\bf Statement: }Then, for orthonormal vectors $N_1$ and $N_2$ spanning the normal space at the point $X(0,0),$ there exists a constant $\Theta=\Theta(h_0r,h_1r^{1+\alpha},h_2r,d_0,\sin\omega,\alpha)\in(0,+\infty)$ such that holds
  $$\kappa_{\Sigma,1}(0,0)^2+\kappa_{\Sigma,2}(0,0)^2
    \le\frac{1}{r^2}\Big\{(h_0r)^2+\Theta\Big\}$$
for the principle curvatures $\kappa_{\Sigma,1}$ and $\kappa_{\Sigma,2}$ in direction $N_\Sigma.$
\end{theorem}
{\it Proof: }The proof follows \cite{Sauvigny_01}, \S\,2. It holds
  $$\begin{array}{lll}
      \kappa_{\Sigma,1}(0,0)^2+\kappa_{\Sigma,2}(0,0)^2\negthickspace
      & = & \negthickspace\displaystyle
            4H_\Sigma(0,0)^2-2K_\Sigma(0,0)
            \,\le\,h_0^2+2|K_\Sigma(0,0)| \\[0.2cm]
      & = & \negthickspace\displaystyle
            \frac{1}{r^2}\,\Big\{(h_0r)^2+r^2|K_\Sigma(0,0)|\Big\}
    \end{array}$$
for $\Sigma=1,2.$ The desired estimate follows from
  $$K_\Sigma(0,0)
    =\frac{(X_{uu}\cdot N_\Sigma)(X_{vv}\cdot N_\Sigma)-(X_{uv}\cdot N_\Sigma)^2}{W^2}\,\Big|_{(0,0)}$$
w.r.t. a normal basis $\{N_1,N_2\}.$ Thus, our proof consists of two steps: First, we have to find a lower bound for the surface area element, and then an upper bound for the second derivatives of the immersion.
\begin{itemize}
\item[1.]
In the first part we will prove the estimate
$$\frac{W(w)}{r^2}\ge C_1 \quad \mbox{for}\; w\in B_{\frac{1}{2}}(0,0)$$
with a constant $C_1=C_1(h_0r,d_0,\sin\omega)>0$.
\begin{itemize}
\item[1.1]
We fix a normal basis of the surface. Note that the vectors $e_3:=(0,0,1,0)$ and $e_4:=(0,0,0,1)$ are not in any tangent plane of the surface due to the fact that
  $$\frac{1}{\sqrt{1+|\nabla\varphi|^2}}\,(-\varphi_x,-\varphi_y,1,0),\quad
    \frac{1}{\sqrt{1+|\nabla\psi|^2}}\,(-\psi_x,-\psi_y,0,1)$$
are normal to the surface (in the corresponding surface point). We see that the inner product of the first vector by $e_3$ and the inner product of the second vector by $e_4$ do not vanish. Therefore, the projections
  $$N_1^*:=e_3-\frac{e_3\cdot X_u^t}{|X_u|^2}\,X_u-\frac{e_3\cdot X_v^t}{|X_v|^2}\,X_v\,,\quad
    N_2^*:=e_4-\frac{e_4\cdot X_u^t}{|X_u|^2}\,X_u-\frac{e_4\cdot X_v^t}{|X_v|^2}\,X_v$$
are a basis of the normal space and can be transformed into an 
orthonormal basis $\{N_1,N_2\}$ of the normal space.
\item[1.2]
In the first part of the proof we will work with the above orthonormal basis $\{N_1,N_2\}.$  Using the differential equation (see Remark 1 above)
  $$\triangle X=2{\mathcal H}(X,N_1)WN_1+2{\mathcal H}(X,N_2)WN_2\quad\mbox{in}\ B,$$
Assumption (A1) yields the estimate
  $$|\triangle X(u,v)|\le 2h_0|\nabla X(u,v)|^2\quad\mbox{in}\ B.$$
The special structure of this differential inequality - the quadratic growth in the gradient on the right hand side - allows us to apply the methods of \cite{Heinz_01}.\\[0.1cm]
We cite two important consequences of our assumptions.
\item[1.3]
From Assumption (A3) we conclude:\\[0.1cm]
Let $\Gamma(B)$ denote the set of all continuous and piecewise differentiable curves $\gamma:[0,1]\to\overline B,$ such that $\gamma(0)=(0,0)$ and $\gamma(1)\in\partial B.$ Then, it holds
  $$\inf_{\gamma\in\Gamma(B)}\int\limits_0^1\left|\frac{d}{dt}\,X\circ\gamma(t)\right|\,dt\ge r.$$
For a proof we refer the reader to \cite{Sauvigny_01}.
\item[1.4]
Assumption (A5) enables us to interpolate the surface area element $W=W(u,v)$ in terms of $|\nabla x^1|^2:$ Namely, it holds
  $$|\nabla x^1|^2\ge W\sin^2\omega\quad\mbox{in}\ B.$$
The proof can be extracted from \cite{Osserman_01}, Lemma 1.1 (note that this result of Osserman makes only use of the conformal representation of the surface).\\[0.1cm]
Now, we estimate the surface area element. For this, we define several two-dimensional auxiliary functions and apply Heinz' results for elliptic systems in the plane (cp. \cite{Heinz_01}). 
\item[1.5]
Let $F^*(u,v):=(x^1(u,v),x^2(u,v))\,:\,\overline B\rightarrow\mathbb R^2$ denote the plane mapping w.r.t. $X=X(u,v).$ Then, for all $w\in B$ there hold
\begin{itemize}
\item[(i)]
$\displaystyle
 |\triangle F^*(w)|
 \le\frac{4h_0}{\sin^2\omega}\,|\nabla F^*(w)|^2$
\item[(ii)]
$\displaystyle
 |\nabla X(w)|^2
 \le\frac{2}{\sin^2\omega}\,|\nabla F^*(w)|^2\,.$
\end{itemize}
For the proof we note that the second inequality follows from $W+W=|\nabla X|^2$ and from the Osserman inequality in 1.4. The differential inequality 
$|\triangle X|\le 2h_0|\nabla X|^2=4h_0W$ (see 1.2) yields
  $$|\triangle F^*|
    \le|\triangle X|
    \le 4h_0W
    \le\frac{4h_0}{\sin^2\omega}\,|\nabla x^1|^2
    \le\frac{4h_0}{\sin^2\omega}\,|\nabla F^*|^2\,.$$
\item[1.6]
Let $w_0\in B$ and $\nu\in(0,1)$ be given such that $B_{2\nu}(w_0):=\{w\in B\,:\,|w-w_0|<2\nu\}\subset B.$ We consider the mapping
  $$Y(w):=\frac{1}{r}\,\{X(w_0+2\nu w)-X(w_0)\},\quad w\in\overline B,$$
as well as the corresponding plane mapping $F(w):=(y^1(w),y^2(w))\,:\,\overline B\rightarrow\mathbb R^2\,.$ The immersion $Y$ satisfies $|Y_u(w)|^2=\frac{4\nu^2}{r^2}\,W(w_0+2\nu w)=|Y_v(w)|^2,$ $Y_u(w)\cdot Y_v(w)^t=0$ and
  $$|\triangle Y(w)|\le2(h_0r)|\nabla Y(w)|^2\quad\mbox{in}\ B.$$
\item[1.7]
From inequality 1.5 and 1.6 we infer
  $$\begin{array}{lll}
      |\triangle F(w)|\negthickspace
      & \le & \negthickspace\displaystyle
              |\triangle Y(w)|
              \,\le\,2(h_0r)|\nabla Y(w)|^2 \\[0.3cm]
      &  =  & \negthickspace\displaystyle
              \frac{16\nu^2(h_0r)}{r^2}\,W(w_0+2\nu w)
              \,=\,\frac{8\nu^2(h_0r)}{r^2}\,|\nabla X(w_0+2\nu w)|^2 \\[0.6cm]
      & \le & \negthickspace\displaystyle
              \frac{16\nu^2(h_0r)}{r^2\sin^2\omega}\,|\nabla F^*(w_0+2\nu w)|^2
              \,=\,\frac{16\nu^2(h_0r)}{r^2\sin^2\omega}\,\frac{r^2}{4\nu^2}\,|\nabla F(w)|^2 \\[0.6cm]
      & \le & \negthickspace\displaystyle
              \frac{4(h_0r)}{\sin^2\omega}\,|\nabla F(w)|^2
    \end{array}$$
for all $w\in B.$ Furthermore, by the assumption that the immersion represents a graph of the form $(x,y,\varphi(x,y),\psi(x,y)),$ $(x,y)\in\overline B_R,$ we conclude that the plane mapping $F=F(u,v)$ is one-to-one and has positive Jacobian $J_F(w)>0$ for all $w\in B$ (see e.g. \cite{Schulz_01}, \cite{Sauvigny_02}). Additionally, Assumption (A4) gives
  $${\mathcal D}[F]\le{\mathcal D}[Y]\le\frac{1}{r^2}\,{\mathcal D}[X]\le d_0$$
for the Dirichlet integrals. We now apply \cite{Heinz_01}, Theorem 6, page 254, which gives the following inner gradient estimate: There is a constant $c_1=c_1(h_0r,d_0,\sin\omega)\in(0,+\infty)$ such that
  $$|\nabla F(u,v)|\le c_1(h_0r,d_0,\sin\omega)\quad\mbox{for all}\ (u,v)\in B_\frac{1}{2}(0,0).$$
\item[1.8]
From 1.7 we infer $\frac{4\nu^2}{r^2}\,W(w_0+2\nu w)\le\frac{1}{\sin^2\omega}\,|\nabla F(w)|^2$ for $w\in B.$ Thus, we arrive at the following estimate
  $$\frac{1}{r^2}\,W(w)
    \le\frac{1}{4\nu^2\sin^2\omega}\,c_1(h_0r,d_0,\sin\omega)
    =:c_2(h_0r,d_0,\sin\omega,\nu)$$
for all $w\in B_{\nu}(w_0)$ and for all $w_0\in B$ such that $B_{2\nu}(w_0)\subset B.$
\goodbreak\noindent
For $w_0\in B_\frac{1}{2}(0,0)$ and $\nu\le\frac{1}{4}$ we get
  $$\frac{1}{r^2}\,W(w)\leq c_2(h_0r,d_0,\sin\omega,\nu)
    \quad\mbox{for all}\ w\in B_{\frac{1}{2}}(0,0).$$
This estimate will be used in the second part of the proof.
\item[1.9]
The properties $J_F(w)>0$ in $B$ and ${\mathcal D}[F]\le d_0,$ as well as the structure of the differential inequality for $|\triangle F(w)|$ (see 1.7) make \cite{Heinz_01}, Lemma 17, page 255 applicable: There exists a constant $c_3=c_3(h_0r,d_0,\sin\omega)\in(0,+\infty)$ such that
  $$|\nabla F(w)|^2\le c_3(h_0r,d_0,\sin\omega)|\nabla F(0,0)|^\frac{2}{5}
    \quad\mbox{for all}\ w\in B_\frac{1}{2}(0,0).$$
It follows, that (cp. 1.6)
  $$\hspace{-0.4cm}
    \begin{array}{lll}
      \displaystyle
      \frac{4\nu^2}{r^2}\,W(w_0+2\nu w)\negthickspace
      & \le & \displaystyle\negthickspace
              \frac{1}{\sin^2\omega}\,|\nabla F(w)|^2
              \,\le\,\frac{c_3(h_0r,d_0,\sin\omega)}{\sin^2\omega}\,|\nabla F(0,0)|^\frac{2}{5} \\[0.6cm]
      & \le & \displaystyle\negthickspace
              \frac{c_3(h_0r,d_0,\sin\omega)}{\sin^2\omega}\,|\nabla Y(0,0)|^\frac{2}{5}
              \,=\,\frac{c_3(h_0r,d_0,\sin\omega)}{\sin^2\omega}
                   \left[\frac{8\nu^2}{r^2}\,W(w_0)\right]^\frac{1}{5}\,.
    \end{array}$$
Rearranging yields the important inequality of Harnack type
  $$\left[\frac{W(w_0)}{r^2}\right]^\frac{1}{5}
    \ge\frac{8^\frac{17}{15}\nu^\frac{8}{5}\sin^2\omega}{c_3(h_0r,d_0,\sin\omega)}\,\frac{W(w_0+2\nu w)}{r^2}
    \quad\mbox{for all}\ w\in B_\frac{1}{2}(0,0),$$
or, which is equivalent,
  $$c_4(h_0r,d_0,\sin\omega,\nu)\left[\frac{W(w)}{r^2}\right]^5\le\frac{W(w_0)}{r^2}
    \quad\mbox{for all}\ w\in B_\nu(w_0)$$
with the constant $c_4(h_0r,d_0,\sin\omega,\nu):=\frac{8^\frac{17}{3}\nu^8\sin^{10}\omega}{c_3(h_0r,d_0,\sin\omega)^5}\in(0,+\infty).$
\end{itemize}
\item[1.10]
Assumption (A4) also ensures that we can estimate the surface area element 
in at least one point: There exists a $w^*\in B_{1-\nu_0}(0,0),$ $\nu_0:=\min(e^{-4\pi d_0},\frac{1}{2})$ such that
  $$\frac{W(w^*)}{r^2}\ge\frac{1}{4(1-e^{-4\pi d_0})}=:c_5(d_0)>0 \, .$$
The constant $\nu_0$ arises from an application of the Courant-Lebesgue lemma (see \cite{Sauvigny_01}).
\item[1.11]
We now show an estimate of the surface area element:\\[0.1cm]
We set $\nu:=\frac{1}{2}\nu_0\in(0,\frac{1}{4}]$ and choose an integer $n=n(\nu)\in\mathbb N$ such that 
$1-2\nu\le \frac{n}{2}\nu\le 1-\nu.$ 
For an arbitrary $w_0\in B_{1-\nu_0}(0,0)$ we define the following points
$$ w_j:=\frac{j}{n}w^*+\frac{n-j}{n}w_0 \quad \mbox{for}\; j=0,\ldots,n \, .$$
Together with 1.10 we have $|w_j|\leq\frac{j}{n}\,|w^*|+\frac{n-j}{n}|w_0|<1-\nu_0,$ and therefore $ B_{2\nu}(w_j)= B_{\nu_0}(w_j)\subset B \, .$ Furthermore, we have
  $$|w_{j+1}-w_j|=|\frac{1}{n}w^*-\frac{1}{n}w_0|\leq\frac{1}{n}|w^*-w_0|
\leq\frac{2(1-\nu_0)}{n}\leq \nu.$$
This implies $ w_{j+1}\in B_{\nu}(w_j) \quad \mbox{for} \; j=0,\ldots,n-1 \, .$
\goodbreak\noindent
We apply the Harnack inequality from 1.9 and obtain
$$\frac{W(w_0)}{r^2}
    \ge c_4\left[\frac{W(w_1)}{r^2}\right]^5
    \ge c_4^{1+5}\left[\frac{W(w_2)}{r^2}\right]^{5^2}
    \ge\ldots
    \ge c_4^{1+5+5^2+\ldots+5^{n-1}}\left[\frac{W(w_n)}{r^2}\right]^{5^n}\,.$$
Recalling $w_n=w^*$, 1.10 gives
  $$\frac{W(w_0)}{r^2}
    \ge c_4^{1+5+5^2+\ldots+5^{n-1}}c_5(d_0)^{5^n}
    =:C_1(h_0r,d_0,\sin\omega)>0$$
for all $w_0\in B_{1-\nu_0}(0,0)$. From $\nu_0\leq\frac{1}{2}$ we now conclude
$$\frac{W(w)}{r^2}\geq C_1(h_0r,d_0,\sin\omega) \quad \mbox{for all}\; 
w\in B_{\frac{1}{2}}(0,0)\, .$$
This ends the first part of the proof.
\item[2.]
In the second part of the proof we have to estimate the second derivatives
of $X$ using the differential equation
$$\triangle X=2{\mathcal H}(X,N_1)WN_1+2{\mathcal H}(X,N_2)WN_2 \, .$$
To do so, we have to give H\"older estimates of the right hand side
of this equation. Especially, we have to construct an orthonormal
basis $\{N_1,N_2\}$ of the normal space at each point, whose
H\"older norm can be estimated. 
\item[2.1]
We define the auxiliary function
$$Z(u,v)=\frac{1}{r}\,\{X(u,v)-X(0,0)\}=\frac{1}{r}\,X(u,v),\quad(u,v)\in\overline B \, .$$
Denoting by $W_Z$ the surface area element of $Z,$ we have $|Z_u|^2=W_Z=|Z_v|^2$ and $Z_u\cdot Z_v^t=0$ in $B.$ It holds $r^2W_Z=W_X$ with $W_X:=|X_u|^2=|X_v|^2.$ We calculate
  $$\begin{array}{lll}
      \triangle Z\negthickspace
      & = & \negthickspace\displaystyle
            \frac{2}{r}\,{\mathcal H}(X,N_1)W_XN_1+\frac{2}{r}\,{\mathcal H}(X,N_2)W_XN_2 \\[0.4cm]
      & = & \negthickspace\displaystyle
            2r{\mathcal H}(rZ,N_1)W_ZN_1+2r{\mathcal H}(rZ,N_2)W_ZN_2\,.
    \end{array}$$
\item[2.2]
Due to 1.8 we have the estimate
  $$|\triangle Z(w)|\le 4(rh_0)c_2(h_0r,d_0,\sin\omega)
    \quad\mbox{for all}\ w\in B_{\frac{1}{2}}(0,0).$$
Furthermore, we get
  $$|Z(u,v)|
    =|Z(u,v)-Z(0,0)|
    \le\frac{1}{2}\,|\nabla Z(\widetilde w)|
    \le\sqrt{c_2(h_0r,d_0,\sin\omega)}
    \quad\mbox{in}\ B_{\frac{1}{2}}(0,0).$$
Now, by potential theoretic estimates, there exists a constant
$c_6(h_0r,d_0,\sin\omega,\alpha)$ such that the H\"older estimate
  $$|Z_{u^i}(w_1)-Z_{u^i}(w_2)|\le c_6(h_0r,d_0,\sin\omega,\alpha)|w_1-w_2|^\alpha
    \quad\mbox{for}\ w_1,w_2\in B_{\frac{1}{4}}(0,0),
    \quad i=1,2$$
holds true. Therefore, for the surface area element it holds
  $$|W_Z(w_1)-W_Z(w_2)|\le c_7(h_0r,d_0,\sin\omega,\alpha)|w_1-w_2|^\alpha
    \quad\mbox{for all}\ w_1,w_2\in B_{\frac{1}{4}}(0,0)$$
with the constant $c_7:=4\sqrt{c_2}\,c_6.$
\item[2.3]
Using the mean value theorem we have the Lipschitz estimate
  $$|Z(w_1)-Z(w_2)|\le 2\sqrt{c_2(h_0r,d_0,\sin\alpha)}\,|w_1-w_2|
    \quad\mbox{for}\ w_1,w_2\in B_{\frac{1}{2}}(0,0).$$
In a certain neighborhood of the origin, we now construct an
orthonormal basis $\{N_1,N_2\}$ of the normal space whose 
H\"older norm can be estimated. \\[1ex]
First we choose vectors $\overline N_1,\overline N_2\in\mathbb R^4$ such that $\overline N_i\cdot Z_{u^j}(0,0)^t=0,$ $\overline N_i\cdot\overline N_j^t=\delta_{ij}\,,$ i.e. $\{\overline N_1,\overline N_2\}$ form a orthonormal basis
of the normal space at the point $Z(0,0)$.
Due to the Gram-Schmidt orthonormalization, for $k=1,2$ we define vectors
  $$N_k^*(w)
    :=\overline N_k
      -\frac{\overline N_k\cdot Z_u(w)}{|Z_u(w)|^2}\,Z_u(w)
      -\frac{\overline N_k\cdot Z_v(w)}{|Z_v(w)|^2}\,Z_v(w)
      \quad\mbox{in}\ B \, .$$
These vectors both belong to the normal space at $Z(w)$ but they may not be
linearly independent. At first, we determine a $\nu_1=\nu_1(h_0r,d_0,\sin\omega,\alpha)>0$ such that
$$|N_k^*(w)|^2
    =1-\frac{[\overline N_k\cdot Z_u(w)]^2}{W_Z(w)}-\frac{[\overline N_k\cdot Z_v(w)]^2}{W_Z(w)}
    \ge\frac{1}{2}
    \quad\mbox{in}\ B_{\nu_1}(0,0).$$
This is possible because, first, we calculate
  $$\begin{array}{lll}
      |\overline N_k\cdot Z_{u^\ell}(w)|^2\negthickspace
      &  =  & \negthickspace\displaystyle
              |\overline N_k\cdot\{Z_{u^\ell}(w)-Z_{u^\ell}(0,0)\}|^2
              \,\le\,|Z_{u^\ell}(w)-Z_{u^\ell}(0,0)|^2 \\[0.2cm]
      & \le & \negthickspace\displaystyle
              c_6(h_0r,d_0,\sin\omega,\alpha)^2|w|^{2\alpha}
    \end{array}$$
by 2.2, and next, by 1.11, we have the lower bound
  $$W_Z(w)\ge C_1(h_0r,d_0,\sin\omega)
    \quad\mbox{in}\ B_{\frac{1}{2}}(0,0).$$
\item[2.4]
We remark that the vectors $N_k^*(w),$ $k=1,2,$ are H\"older continuous in $B_{\nu_1}(0,0)$
and the H\"older estimate
$$ |N_k^*(w_1)-N_k^*(w_2)|\leq c_8(h_0r,d_0,\sin\omega,\alpha)|w_1-w_2|^\alpha\,,
\quad \ w_1,w_2\in B_{\nu_1}(0,0),$$
holds true with a constant $c_8(h_0r,d_0,\sin\omega,\alpha)$ coming from the H\"older estimate
for $Z_{u_j}$ and the lower bound of $W_Z.$ Now, for $k=1,2$ we define
$$\widetilde N_k(w):=\frac{N_k^*(w)}{|N_k^*(w)|}\quad\mbox{in}\ B_{\nu_1}(0,0).$$
These vectors are well defined because of $|N_k^*(w)|^2\ge\frac{1}{2}$ in $B_{\nu_1}(0,0),$ 
but they are not orthogonal. Note that
$$N_1^*\cdot N_2^*
    =-\frac{(\overline N_1\cdot Z_u^t)(\overline N_2\cdot Z_u^t)}{W_Z}
     -\frac{(\overline N_1\cdot Z_v^t)(\overline N_2\cdot Z_v^t)}{W_Z}\,,$$
and, therefore, due to 2.3
$$|\widetilde N_1\cdot\widetilde N_2|=\frac{|N_1^*\cdot N_2^*|}{|N_1^*||N_2^*|}
  \le\frac{2}{C_1}\,
   (|\overline N_1\cdot Z_u^t||\overline N_2\cdot Z_u^t|+
  |\overline N_1\cdot Z_v^t||\overline N_2\cdot Z_v^t|)
    \le\frac{4c_6^2}{C_1}\,|w|^{2\alpha}\,.$$
Thus we can find a $\nu_2=\nu_2(h_0r,d_0,\sin\omega,\alpha)$ with
$0<\nu_2\leq\nu_1$ such that
$$ |\widetilde N_1\cdot\widetilde N_2|\leq \frac{1}{2} \quad \mbox{in} \ B_{\nu_2}(0,0) \, .$$
Now we define the vectors
$$N_1(w):=\widetilde N_1(w),\quad
  N_2(w):=\frac{\widetilde N_2(w)-(N_1\cdot\widetilde N_2)N_1}{\sqrt{1-(N_1\cdot\widetilde N_2)^2}}
 \quad \mbox{in}\ B_{\nu_2}(0,0)\, .$$
$N_2$ is well defined because for its denominator we have $ 1-(N_1\cdot \widetilde N_2)^2=1-(\widetilde N_1\cdot\widetilde N_2)^2\geq \frac{3}{4} \, .$ The vectors $N_1(w),N_2(w)$ are an orthonormal basis of the normal space at
each point $Z(w)$ for $w\in B_{\nu_2}(0,0)$. Furthermore, the H\"older estimate
$$ |N_k(w_1)-N_k(w_2)|\leq c_9(h_0r,d_0,\sin\omega,\alpha) |w_1-w_2|^\alpha 
\quad \mbox{for}\ w_1,w_2\in B_{\nu_2}(0,0)$$
holds true with a constant $c_9(h_0r,d_0,\sin\omega,\alpha)$ which can be calculated
directly using the H\"older estimate for $N^*_k$.
\item[2.5]
We now use the differential equation (recall Remark 1 to Definition 2)
$$ \triangle Z=2r\mathcal H(rZ,N_1)W_Z N_1+2r\mathcal H(rZ,N_2)W_Z N_2 
\quad \mbox{in} \; B_{\nu_2}(0,0) \, .$$
We already showed the estimate $ |\triangle Z(w)|\leq 4(h_0r)c_2(h_0r,d_0,\sin\omega)$ in $B_{\nu_2}(0,0).$ Using the assumptions on $\mathcal H$ we obtain the H\"older estimate
  $$\begin{array}{l}
      |\mathcal H(rZ(w_1),N_k(w_1))-\mathcal H(rZ(w_2),N_k(w_2))| \\[0.2cm]
      \hspace{0.6cm}
      \leq\,h_1 r^\alpha|Z(w_1)-Z(w_2)|^\alpha+h_2|N_k(w_1)-N_k(w_2)| \\[0.2cm]
      \hspace{0.6cm}
      \leq\,h_12^\alpha r^\alpha c_2^\frac{\alpha}{2}\,|w_1-w_2|^\alpha+h_2 c_9|w_1-w_2|^\alpha\,.
    \end{array}$$
Thus we can find a constant $c_{10}=c_{10}(h_0r,h_1r^{1+\alpha},h_2r,d_0,\sin\omega,\alpha)$
such that
$$|\triangle Z(w_1)-\triangle Z(w_2)|\leq c_{10} |w_1-w_2|^\alpha
\quad \mbox{for} \; w_1,w_2\in B_{\nu_2}(0,0) \, .$$
Setting $\nu_3:=\frac{1}{2}\nu_2$, the interior Schauder estimates 
give a constant $C_2\in(0,+\infty)$ such that
$$|Z_{uu}(w)|,|Z_{uv}(w)|,|Z_{vv}(w)|\le C_2(h_0r,h_1r^{1+\alpha},h_2r,d_0,\sin\omega,\alpha)\quad\mbox{in}\ B_{\nu_3}(0,0).$$
From the beginning of the proof we recall
  $$\kappa_{1,\Sigma}(0,0)^2+\kappa_{2,\Sigma}(0,0)^2
    \le\frac{1}{r^2}\,\left\{(h_0r)^2+\frac{|Z_{uu}(0,0)||Z_{vv}(0,0)|+|Z_{uv}(0,0)|^2}{W_Z(0,0)^2}\right\}.$$
Setting $\Theta(h_0r,h_1r^{1+\alpha},h_2r,d_0,\sin\omega,\alpha):=\frac{2C_2^2}{C_1^2}$ we arrive at
  $$\kappa_{1,\Sigma}(0,0)^2+\kappa_{2,\Sigma}(0,0)^2
    \le\frac{1}{r^2}\Big\{(h_0r)^2+\Theta\Big\}.$$
This proves the statement.\hfill{$\Box$}
\end{itemize}
\vspace{1cm}
{\it Remarks:}
\begin{itemize}
\item[1.]
For the introduction of conformal parameters into the graph we refer to \cite{Schulz_01} and \cite{Sauvigny_02}. In the theorem, we can replace the domain $B_R$ of the graph  by a $C^2$-Jordan domain $\Omega$ using Riemann's mapping theorem (see e.g. \cite{Sauvigny_03}).
\item[2.]
Assumption $(A5)$ is basically needed to establish the inequality
$$ |\triangle F^*(w)|\leq \frac{4h_0}{\sin^2\omega}\,|\nabla F^*(w)|^2 $$
for the plane mapping $F^*$. This is the difference to estimates for
surfaces in $\mathbb R^3,$ where $ |\triangle F^*(w)|\leq 4h_0\,|\nabla F^*(w)|^2 $ holds due to the conformal representation (note, that in this case $|X^3_u|^2\leq |X^1_v|^2+|X^2_v|^2$).
\item[3.]
The proof does not make use of the special dimension $n=4$.
Thus, it can be carried over to the more general case of dimension $n\geq 4$.
\end{itemize}
In the case of minimal surfaces, i.e. $\mathcal H\equiv 0$, the proof yields a result of Bernstein type since the a priori constant $\Theta$ does not depend on $r$ any more. For $r\to\infty$ we obtain the
\begin{corollary}
A complete minimal graph $X(x,y)=(x,y,\varphi(x,y),\psi(x,y))$ which satisfies assumptions (A4) and (A5) is an affine plane.
\end{corollary}
Our results are motivated from \cite{Osserman_01}, where curvature estimates and a Bernstein-type result for minimal immersions in $\mathbb R^n$ are proved.
\vspace*{0.8cm}
\noindent
Matthias Bergner, Steffen Fr\"ohlich\\
Technische Universit\"at Darmstadt\\
Fachbereich Mathematik, AG 4\\
Schlo\ss{}gartenstra\ss{}e 7\\
D-64289 Darmstadt\\
Germany\\[0.2cm]
e-mail: bergner@mathematik.tu-darmstadt.de\\
\hspace*{1.12cm}
sfroehlich@mathematik.tu-darmstadt.de
\end{document}